\documentclass[reqno]{amsart}
\usepackage{latexsym}
\usepackage{stmaryrd}
\usepackage{epsfig}
\theoremstyle{plain}
\usepackage{txfonts}

\def\XXint#1#2#3{{\setbox0=\hbox{$#1{#2#3}{\int}$}
     \vcenter{\hbox{$#2#3$}}\kern-.5\wd0}}

\makeglossary
\newcommand{\ca}[1]{Card\lt(#1\rt)}

\newcommand{\ep}{\epsilon}

\newcommand{\lt}{\left}
\newcommand{\lm}{\lambda}

\newcommand{\cv}{\mathcal{C}}
\newcommand{\la}{\langle}
\newcommand{\ra}{\rangle}
\newcommand{\rt}{\right}
\newcommand{\nl}{\newline}

\newcommand{\nn}{\nonumber}

\newcommand{\ti}{\tilde}

\newcommand{\qd}{\quad}

\newcommand{\TT}{\mathcal{T}}

\newcommand{\spt}{\mathrm{Spt}}

\newcommand{\veps}{\varphi}
\newcommand{\vva}{\varrho}
\newcommand{\vthe}{\vartheta}
\newcommand{\ui}{\tilde{u}}
\newcommand{\wt}{\widetilde}

\newcommand{\R}{\mathrm {I\!R}}

\newcommand{\divv}{\mathrm{div}}
\newtheorem{a1}{Lemma}
\newtheorem{a2}{Theorem}

\begin{document}
\title[A simple proof of the characterization of functions of low Aviles Giga energy]
{A simple proof of the characterization of functions of low Aviles Giga energy on a ball via regularity}  
\author{Andrew Lorent}
\address{Mathematics Department\\University of Cincinnati\\ 2600 Clifton Ave.\\ Cincinnati\\ Ohio 45221\\}
\email{lorentaw@uc.edu}\date{\today}
\maketitle

\begin{abstract}

The Aviles Giga functional is a well known second order functional that forms a model for blistering and in a certain 
regime liquid crystals, a related functional models thin magnetized films. Given 
Lipschitz domain $\Omega\subset \R^2$ the 
functional is $I_{\ep}(u)=\frac{1}{2}\int_{\Omega} \ep^{-1}\lt|1-\lt|Du\rt|^2\rt|^2+\ep\lt|D^2 u\rt|^2 dz$ where 
$u$ belongs to the subset of functions
in $W^{2,2}_{0}(\Omega)$ whose gradient (in the sense of trace) satisfies $Du(x)\cdot \eta_x=1$ where 
$\eta_x$ is the inward pointing unit normal to $\partial \Omega$ at $x$.

In \cite{otto} Jabin, Otto, Perthame characterized a class of functions which includes all limits of sequences $u_n\in W^{2,2}_0(\Omega)$ 
with $I_{\ep_n}(u_n)\rightarrow 0$
as $\ep_n\rightarrow 0$. A corollary to their work is that if there exists such a sequence $(u_n)$ for a bounded
domain $\Omega$, then $\Omega$ must be a ball and (up to change of sign) $u:=\lim_{n\rightarrow \infty} u_n
=\mathrm{dist}(\cdot,\partial\Omega)$. Recently \cite{lor12} we provided a quantitative generalization of this corollary over the 
space of convex domains using `compensated compactness' inspired calculations of \cite{mul2}. In this note we 
use methods of regularity theory and ODE to provide a sharper estimate and a much simpler proof for the case where $\Omega=B_1(0)$ without the requiring the trace condition on $Du$. 
\end{abstract}

\section{Introduction}
Let
\begin{equation}
\label{aaz1}
I_{\ep}\lt(u\rt):=\int_{\Omega}
\ep^{-1}\lt|1-\lt|Du\rt|^2\rt|^2+\ep\lt|D^2 u\rt|^2 dz. 
\end{equation}

The functional $I_{\ep}$ forms a model for blistering and (in
certain regimes) for a model for liquid crystals \cite{avgig}, \cite{kohn}. 
In addition there is a closely related functional modeling thin magnetic
films \cite{mul2}, \cite{deott1}, \cite{cdeott}, \cite{ser1}, \cite{ser2}. For function $u\in W_0^{2,2}(\Omega)$ we 
refer to $I_{\ep}(u)$ as the \em Aviles Giga energy \rm of $u$. 

For an example of a candidate minimizer take the distance function 
from the boundary $\psi(x):=\mathrm{dist}(x,\partial \Omega)$ convolved by a standard convolution kernel $\rho_{\ep}$ with support 
of diameter $\ep$. It has been conjectured that for convex domains $\Omega$, the 
minimizers of $I_{\ep}$ have the structure suggested by this construction, i.e.\ they are 
in some quantitative sense close to the distance function from the boundary, Section 5.3 \cite{gioort},\cite{avgig0}.    

The first progress on this conjecture was achieved by Jin, Kohn \cite{kohn} whose showed that if $I_{\ep}$ is minimized over 
\begin{equation}
\label{aab12}
\Lambda\lt(\Omega\rt):=\lt\{\begin{array}{ll}
v\in W^{2,2}_{0}\lt(\Omega\rt):&\frac{\partial v}{\partial \eta_z}=1\text{ where }\eta_z\text{ is the inwards }\\
&\text{pointing unit normal to }\partial \Omega\text{ at }z\end{array}\rt\}
\end{equation}
where $\Omega$ is taken to be an ellipse then as $\ep\rightarrow 0$ the energy of the minimizer of $I_{\ep}$ tends to the 
energy of $\psi*\rho_{\ep}$. Their method was to take arbitrary 
$u\in \Lambda(\Omega)$ and to construct vectors 
fields $\Sigma_1$, $\Sigma_2$ out of third order polynomials of the partial derivatives of $u$ that have the property that the 
divergence of these vectors fields is bounded above by $I_{\ep}(u)$. Using the trace condition $\frac{\partial u}{\partial \eta}=1$ and 
the fact that $\Omega$ is an ellipse the lower bound provided by the  
divergence of $\Sigma_1,\Sigma_2$ can be explicitly calculated and shown to be asymptotically sharp as $\ep\rightarrow 0$. 

As has been discussed in \cite{kohn}, \cite{avgig0}, \cite{amb2.5} the functional $I_{\ep}$ minimized over $W^{2,2}_0(\Omega)$ has many 
features in common with the functional $J^p(v)=\int_{J_{Dv}} \lt|Dv^{+}-Dv^{-}\rt|^p dH^1$ for the case $p=3$, when minimized over the space 
$Dv\in BV(\Omega)$ with $\lt|Dv(x)\rt|=1$ a.e.\ $x$ and $v=0$ on $\partial \Omega$. Aviles Giga \cite{avgig1} showed that if $\Omega$ is convex and 
polygonal then the distance function is the minimizer of $J^1$ over the subspace of piecewise affine functions satisfying 
these conditions. They conjectured the same is true for $p=3$.

From a somewhat different direction a strong result has been proved \cite{otto} by Jabin, Otto, Perthame  who characterized a class of 
functions which includes all limits of sequences $u_n\in W^{2,2}_0(\Omega)$ with $I_{\ep_n}(u_n)\rightarrow 0$
as $\ep_n\rightarrow 0$. A corollary to their work is that if there exists such a sequence $(u_n)$ for a bounded
domain $\Omega$, then $\Omega$ must be a ball and (up to change of sign) $u:=\lim_{n\rightarrow \infty} u_n
=\mathrm{dist}(\cdot,\partial\Omega)$. In \cite{lor12}, a quantitative generalization of this corollary was achieved for 
the class of bounded convex domains, a corollary to the main result of \cite{lor12} is the following.  
\begin{a2}[Lorent 2009]
\label{CC1}
Let $\Omega$ be a convex set with diameter $2$, $C^2$ boundary and curvature bounded above by $\ep^{-\frac{1}{2}}$. Let $\Lambda(\Omega)$ be defined 
by (\ref{aab12}). There exists 
positive constants $C>1$ and $\lambda<1$ such that if $u$ is a minimizer of $I_{\ep}$ over $\Lambda(\Omega)$, then 
\begin{equation}
\label{uz9}
\|u-\zeta\|_{W^{1,2}(\Omega)}\leq C \lt(\ep+\inf_{y}\lt|\Omega\triangle B_1(y)\rt|\rt)^{\lambda}
\end{equation}
where $\zeta(z)=\mathrm{dist}(z,\partial \Omega)$. 
\end{a2}
We take constant $\lambda=\frac{1}{2731}$ and thus the control represented by 
inequality (\ref{uz9}) is far from optimal.  Theorem \ref{CC1} follows from 
Theorem 1 of \cite{lor12} which is a characterization of domains $\Omega$ and functions $u$ for which the Aviles Energy is small, more specifically 
there exists a constant $\gamma$ such that given $u\in \Lambda(\Omega)$ such that $I_{\ep}(u)=\beta$ then 
$\lt|\Omega\triangle B_1(0)\rt|\leq c\beta^{\gamma}$ and $\int_{B_1(0)} \lt|Du(z)+\frac{z}{\lt|z\rt|}\rt|^2 dz\leq c\beta^{\gamma}$, 
here we can take $\gamma=512^{-1}$. The proof of Theorem 1 of \cite{lor12} is fairly involved, it relies heavily on the characterization 
of `entropies' for the Aviles Giga energy that was achieved in \cite{mul2}, (see Lemma 3). While 
the calculations in \cite{lor12} are elementary and self contained, they can appear quite unmotivated to those unfamiliar with the background of 
\cite{mul2}. In addition the trace condition on the gradient in the definition of $\Lambda(\Omega)$ is used in an essential way. 

The proof of Theorem 1 requires quite a careful construction of an upper bound of the Aviles Giga energy of a minimizer on a domain with smooth boundary 
that is `close' to a ball, then the theorem follows by application of Theorem 1 \cite{lor12}. The many steps required to complete the proof result in a 
gradual loss of control resulting in the constant $\lm=\frac{1}{2731}$.  

The propose of this note is twofold, firstly 
to provide a simple proof of a characterization of the minimizers of the Aviles Giga energy on a ball with a sharper estimate 
and secondly to prove the result without the trace condition on the gradient, specifically to characterize the minimizers 
over $W^{2,2}_0(B_1(0))$.  Additionally we find it worthwhile to introduce new methods to study the characterization of minimizers of 
$I_{\ep}$, the regularity theory and ODE approach of this note is quite different from previous methods of   \cite{avgig1}, \cite{kohn}, \cite{otto}, \cite{lor12}. 
Our main theorem is; 

\begin{a2} 
\label{T2}
Let $u$ be a minimizer of $I_{\ep}$ over $W^{2,2}_0(B_1(0))$. Then there exists $\xi\in \lt\{1,-1\rt\}$
$$
\int_{B_1(0)} \lt|Du(x)+\xi\frac{x}{\lt|x\rt|}\rt|^2 dx\leq c\ep^{\frac{1}{6}}(\log(\ep^{-1}))^{\frac{13}{6}}.
$$
\end{a2}

The desirability of a simpler proof with a better estimate has already been discussed, it is of interest to prove a characterization without a 
trace condition on the gradient due to the fact this is a strong assumption that is inappropriate for a number of physical models. More specifically the condition $Du(x)\cdot \eta=1$ for $x\in \partial \Omega$ is not 
natural in the context of blistering, Gioia Ortiz \cite{gioort} proposed instead $Du(x)\cdot \eta_x=0$. The original functional proposed by Aviles 
Giga \cite{avgig0} to study liquid crystals also has this trace condition.  In addition for the micro-magnetic analogue 
of functional $I_{\ep}$ there is nothing like a pointwise condition on the trace, \cite{deott1}, \cite{cdeott}. This micro-magnetic 
functional is given by $M_{\ep}(v)=\ep^{-1}\int_{\R^2} \lt|H(\ti{v})\rt|^2+\ep\int_{\Omega} \lt|Dv\rt|^2$ where $H$ is the Hodge projection onto 
curl free vector fields and $\ti{v}$ is the extension of $v$ to $0$ outside $\Omega$, this functional is minimized over $W^{1,2}(\Omega:S^1)$. As mentioned, in the proof of Theorem 1 \cite{lor12} the trace condition is used in an essential way, this is also true of the proof of Theorem 5.1 
\cite{kohn}. In order to achieve a characterization for less rigid functionals, methods need to be developed that do not use this trace condition. A related but different micro-magnetic functional $E_{\ep}$ was 
studied by Ignat, Otto \cite{igotto}. They also achieved a characterization of minimizers $E_{\ep}$ 
showing that minimizers converge to Neel Walls, the focus of $E_{\ep}$ was to provide a two 
dimensional approximation of the micro-magnetic energy in the absence of an external field and 
crystal anisotropy.

The proof of Theorem \ref{T2} requires establishing the essentially folklore fact that critical points of the Aviles Giga energy have $W^{2,3}$ 
regularity and their gradients satisfy certain natural Caccioppoli inequalities. The much more subtle question of regularity of critical points 
of functional $M_{\ep}$ has been studied by Carbou \cite{carbou} and Hardt, Kinderlehrer \cite{hakin}. The non-local term in $M_{\ep}$ makes the 
Euler Lagrange equation harder 
to study and in some sense weaker regularity has been proved, it is not clear if the Caccioppoli inequalities needed for the proof presented in 
this note are available via the methods of \cite{carbou}. Working with a three dimensional model, 
different methods are used in \cite{hakin} and Caccioppoli inequalities are 
established off a discrete set\footnote{ It appears possible that the methods of \cite{hakin} would establish the appropriate Caccioppoli inequalities 
everywhere in the interior if the arguments were carried through for the two dimensional model, if this is the case the strategy of this note 
would likely yield a characterization of minimizers of $M_{\ep}$ for where $\Omega=B_1(0)$.}.  

Roughly speaking the main open problems related to the Aviles Giga functional are either; (A) conjectures on how the energy concentrates, 
specifically the $\Gamma$-convergence conjecture of \cite{amb2.5} and related problems. Or (B) conjectures about the minimizer of $I_{\ep}$. 
It is know from \cite{kohn} that for non-convex domains the minimizer does not need to be the distance function from the boundary 
(contrast this with the main theorem of \cite{amleri} which showed that for a sequence $\ep_n\rightarrow 0$, 
the minimizer $m_n$ of the micro-magnetics functional $M_{\ep_n}$ must converge to the 
rotated gradient of distance function for any connected open Lipschitz domain). However as mentioned for general convex domains the conjecture remains 
largely open, in \cite{lor12} we developed methods that prove the conjecture for convex domains with low Aviles Giga energy, it is likely these 
methods could be used to prove the same result for general low energy domains with $C^2$ boundary. For domains with Aviles Giga energy of order 
$O(1)$ neither the methods of \cite{lor12} or this note yield much. A very attractive open problem is to characterize the minimizers in the case where 
$\Omega$ is an ellipse, given the sharp lower bound provided by \cite{kohn} in this case there seems to be much concrete information about this problem - yet it appears to be out of reach of current methods.

\section{Proof sketch}

Beyond the regularity issues mentioned in the introduction the proof reduces to essentially applying an ODE and using the Pythagorean Theorem. In order to sketch the main strategy of the proof we will make a number of 
assumptions that we will later show are not needed. 

We start by assuming for a moment that the cardinality of the set of critical points of $Du$ is $1$, i.e.\ 
\begin{equation}
\label{azu1}
\ca{\lt\{x\in B_1(0):\lt|Du(x)\rt|=0\rt\}}=1.
\end{equation}
In addition let us temporarily assume we have the (in the sense of trace) boundary condition 
\begin{equation}
\label{fa1}
Du(x)=-\frac{x}{\lt|x\rt|} \text{ for }x\in \partial B_1(0).
\end{equation}
So let $z_0\in B_1(0)$ be the point for which $\lt|Du(z_0)\rt|=0$. Take $y_0=-z_0 \R\cap \partial B_1(0)$ and let 
$X(0)=y_0$, $\frac{dX}{dt}(s)=Du(X(s))$. For $z\in \lt\{X(s):s\in [0,t]\rt\}$ let $t_z$ denote 
the tangent to this curve at $z$. Now for any $t>0$
$$
u(X(t))=u(X(t))-u(X(0))=\int_{\lt\{X(s):s\in\lt[0,t\rt]\rt\}} Du(z)\cdot t_z dH^1 z.
$$
If we also assume 
\begin{equation}
\label{azu2}
\lt|Du(z)\rt|\eqsim 1\text{ for }z\in\lt\{X(s):s\in\lt[0,t\rt]\rt\}
\end{equation}
then we could conclude that 
$$
\lt|u(X(t))\rt|\eqsim H^1(X(s):s\in\lt[0,t\rt])\geq \lt|X(t)-X(0)\rt|.
$$
Now by (\ref{fa1}) we know that the path $X(t)$ has to run \em into \rm $B_1(0)$ and can not escape this 
domain, so we must have $X(t)\rightarrow z_0$ as $t\rightarrow \infty$ we have $\lt|u(z_0)\rt|\geq \lt|z_0-X(0)\rt|=\lt|z_0\rt|+1$. 

As will be established later in Lemma \ref{L0.5}, $\inf_{v\in W_0^{2,2}(B_1(0))} I_{\ep}(v)\leq 
c\ep\log(\ep^{-1})$. Hence if $u$ is a mimiser of $I_{\ep}$, 
\begin{equation}
\label{avvu5}
\int_{B_1(0)} \lt|1-\lt|Du\rt|^2\rt|^2 dx\leq c\ep^2\log(\ep^{-1}) 
\end{equation}
so we know $u$ `is close to being' $1$-Lipschitz and thus 
$\lt|u(z_0)\rt|\succapprox 1$, hence $\lt|z_0\rt|\eqsim 0$ and $\lt|u(z_0)\rt|\eqsim 1$. Again since $u$ is close to $1$-Lipschitz, 
\begin{equation}
\label{fineq14}
\lt|u(x)\rt|\eqsim 1\text{ for any }x\in B_{\ep^{\frac{1}{4}}}(0).
\end{equation}

Now for $y\in \partial B_1(0)$ let  
$e^x(y)=\int_{\lt[x,y\rt]} \lt|1-\lt|Du\rt|^2\rt| dH^1$. Let $J_x(z)=\frac{z-x}{\lt|z-x\rt|}$, note that 
$\lt|D J_x(z)\rt|\leq \frac{2}{\lt|x-z\rt|}$, so by the Co-area formula 
\begin{eqnarray}
\int_{\partial B_1(0)} e^x(y) dH^1 y&=&\int_{S^1}\int_{J_x^{-1}(\theta)} \lt|1-\lt|Du(z)\rt|^2\rt| dH^1 z dH^1 \theta\nn\\
&=&\int_{B_1(0)} \lt|1-\lt|Du(z)\rt|^2\rt|\lt|D J_x(z)\rt| dz\nn\\
&\leq& c\int_{B_1(0)} \lt|1-\lt|Du(z)\rt|^2\rt|\lt|z-x\rt|^{-1} dz.\nn
\end{eqnarray}
Now by Fubini and (\ref{avvu5}) we have 
$$
\int_{B_{\ep^{\frac{1}{4}}}(0)} \int_{B_1(0)} \lt|1-\lt|Du(z)\rt|^2\rt|\lt|z-x\rt|^{-1} dz dx\leq 
c\ep^{\frac{5}{4}}\sqrt{\log(\ep^{-1})}
$$
thus we can assume we chose $x\in B_{\ep^{\frac{1}{4}}}(0)$ such that $\int_{\partial B_1(0)} e^x(y) dH^1 y\leq 
c\ep^{\frac{3}{4}}\sqrt{\log(\ep^{-1})}$. Now
\begin{eqnarray}
\label{uzuu4}
\int_{\lt[x,y\rt]} \lt|Du(z)+\frac{y-x}{\lt|y-x\rt|}\rt|^2 dH^1 z&=&
\int_{\lt[x,y\rt]} \lt|Du(z)\rt|^2+2 Du(z)\cdot \frac{y-x}{\lt|y-x\rt|}+1 dH^1 z \nn\\
&\leq& 2\lt|x-y\rt|-2u(x)+e^x(y)\nn\\
&\overset{(\ref{fineq14})}{\lessapprox}& e^x(y). 
\end{eqnarray}
So 
\begin{eqnarray}
\label{uuau11}
\int_{B_1(0)} \frac{\lt|Du(z)+\frac{z-x}{\lt|z-x\rt|}\rt|^2}{\lt|z-x\rt|}  dz&\leq&c\int_{y\in \partial B_1(0)} \int_{\lt[x,y\rt]} \lt|Du(z)+\frac{y-x}{\lt|y-x\rt|}\rt|^2 dH^1 z dH^1 y\nn\\
&\overset{(\ref{uzuu4})}{\leq}&c\int_{y\in \partial B_1(0)} e^x(y)\; dH^1 y\nn\\
&\leq & c\ep^{\frac{3}{4}}\sqrt{\log(\ep^{-1})}.
\end{eqnarray}
As for `most' $z\in B_1(0)$, $\lt|\frac{z}{\lt|z\rt|}-\frac{z-x}{\lt|z-x\rt|}\rt|\leq c\ep^{\frac{1}{8}}$ so we have 
$\int_{B_1(0)} \lt|Du(z)+\frac{z}{\lt|z\rt|}\rt|^2 dz\overset{(\ref{uuau11})}{\leq} c\ep^{\frac{1}{8}}$.

Now the big assumptions we made are (\ref{azu1}), (\ref{azu2}) and to a lesser 
extent (\ref{fa1}). The main work of this note is to find substitutes for 
these assumptions.  

What assumption (\ref{azu1}) provides is the existence of a long integral path of the vector field $Du$ which using 
assumption (\ref{azu2}) we can show is close to a straight line. In order to find such a path, it is sufficient 
to show that the set of critical points of $Du$ are merely low in number, using the energy upper bound and regularity of 
minimizers of $I_{\ep}$ that is what we will be able to do. 

Now if we define $v(z)=u(\ep z)$ then $v$ satisfies $\Delta^2 v+\divv\lt(\lt(1-\lt|Dv\rt|^2\rt)Dv\rt)=0$ which is an 
Elliptic equation with right-hand side bounded in $H^{-1,p}(B_{\ep^{-1}}(0))$ for all $p>1$. Thus it is not hard to believe $Dv$ is Holder 
so if $\lt|Dv(z_0)\rt|=0$  for some $z_0$ then there must be a constant $c_0$ such that 
$\sup\lt\{\lt|Dv(z)\rt|:z\in B_{c_0}(z_0)\rt\}\leq \frac{1}{2}$ so after rescaling we have that for every $z_1$ such that 
$\lt|Du(z_1)\rt|=0$ we have that $\sup\lt\{\lt|Du(z)\rt|:z\in B_{c_1 \ep}(z_0)\rt\}\leq \frac{1}{2}$. Thus by (\ref{avvu5}) we 
have that we can have as most $c\log(\ep^{-1})$ critical points of $I_{\ep}$ that are spaced out by $\ep$. So cutting 
$B_1(0)$ into $N=\lt[ \frac{4 c\pi}{\log(\ep^{-1})}\rt]$ equal angles slices which we denote by 
$T_1,T_2,\dots T_N$ then at least half of them do not have any critical points of $Du$. So if $T_1$ is one of them, taking 
$y_0$ to be the center of the arc $T_1\cap \partial B_1(0)$ the ODE $X(0)=y_0$, $\frac{dX}{dt}(s)=Du(X(s))$ has to run 
until it hits $\partial T_1$. 

Now the second main assumption we made is (\ref{azu2}). Again since for minimizer $u$ we 
know that $I_{\ep}(u)\leq c\ep\log(\ep^{-1})$, so 
$$
\int_{B_1(0)} \lt|1-\lt|Du\rt|^2\rt|\lt|D^2 u\rt| dx\leq c\ep \log(\ep^{-1}).
$$ 
Take $v\in S^1$, for all but $c(\ep \log(\ep))^{\frac{1}{3}}$ lines $L$ 
parallel to $v$ we have that $\int_{L} \lt|1-\lt|Du\rt|^2\rt|\lt|D^2 u\rt| dH^1 x \leq (\ep \log(\ep))^{\frac{2}{3}}$. 
Now on the line $L$ if there is a point $z_1\in L$ with $\lt|1-\lt|Du(z_1)\rt|^2\rt|\geq 5(\ep\log(\ep^{-1}))^{\frac{1}{3}}$ then we 
must be able to find $z_2,z_3$ we have $\inf\lt\{\lt|1-\lt|Du(y)\rt|^2\rt|:y\in\lt[z_2,z_3\rt]\rt\}\geq 4(\ep\log(\ep^{-1}))^{\frac{1}{3}}$ and 
$\lt|1-\lt|Du(z_3)\rt|^2\rt|\geq 5 (\ep\log(\ep^{-1}))^{\frac{1}{3}}$, $\lt|1-\lt|Du(z_2)\rt|^2\rt|\leq 4(\ep\log(\ep^{-1}))^{\frac{1}{3}}$ 
then  
$$
(\ep\log(\ep))^{\frac{2}{3}}\geq \int_{z_2}^{z_3} \lt|1-\lt|Du(y)\rt|^2\rt|\lt|D^2 u(y)\rt| dH^1 y\geq 4(\ep\log(\ep^{-1}))^{\frac{1}{3}}
\int_{z_2}^{z_3} \lt|D^2 u(y)\rt| dH^1 y\geq 4(\ep\log(\ep))^{\frac{2}{3}}
$$
which is a contradiction. Thus for most lines $L$ we know that $\sup\lt\{\lt|1-\lt|Du(z)\rt|^2\rt|:y\in L\cap B_1(0)\rt\}\leq 
5(\ep\log(\ep))^{\frac{1}{3}}$. For vector $w\in \R^2$ define $\la w\ra:=\lt\{\lm w:\lm \in \R\rt\}$ and 
given subspace $V$ let $P_{V}$ denote the orthogonal projection onto $V$. For subset $S\subset \R^n$ let 
$\lt|S\rt|$ denote the Lebesgue $n$-measure of $S$. Now if we run an ODE $X(0)=y_0$, $\frac{dX}{dt}(s)=Du(X(s))$ between $0$ and $t$ then taking 
$v=\frac{X(t)-X(0)}{\lt|X(t)-X(0)\rt|}$ then we have a set $G\subset P_{\la v\ra}(\lt[X(0),X(t)\rt])$ with 
$\lt|P_{\la v\ra}(\lt[X(0),X(t)\rt])\backslash G\rt|\leq c (\ep\log(\ep^{-1}))^{\frac{1}{3}}$ and if 
$z\in \lt\{X(s):s\in \lt[0,t\rt]\rt\}\cap P^{-1}_{\la v\ra}(x)$ for some $x\in G$, 
then $\lt|\lt|Du(z)\rt|^2-1\rt|\leq 5(\ep\log(\ep))^{\frac{1}{3}}$ thus 
the part of the path $\lt\{X(s):s\in\lt[0,t\rt]\rt\}$ that is in the set $P_{\la v\ra}^{-1}(G)$ is such that $\lt|Du(z)\rt|\eqsim 1$. 
So the $H^1$ measure of the set of points $x\in \lt\{X(s):s\in\lt[0,t\rt]\rt\}$ for which we can assume $\lt|Du(x)\rt|\eqsim 1$ is of 
measure as least $\lt|X(0)-X(t)\rt|-c(\ep \log(\ep^{-1}))^{\frac{1}{3}}$ and hence assumption (\ref{azu2}) can in effect be justified. It is worth noting that the idea of following integral curves of the vector field 
given by $D u$ (where $u$ is the limit of a sequence of functions whose Aviles Giga energy tends to zero) was used by \cite{otto} and a similar idea later by \cite{igotto}. 

Finally we also assumed (\ref{fa1}), the only purpose of this assumption was to allow us to run 
an ODE starting from $y_0\in \partial B_1(0)$ without it immediately trying to leave the domain. Recall $y_0$ was the 
point at the center of the arc $\partial T_1\cap \partial B_1(0)$. If instead of starting at this point we started 
at $y_0+c\frac{\eta_{y_0}}{(\log(\ep^{-1}))^2}$ then running the ODE forwards and backwards until both 
ends hit $\partial T_1$, then we will have a path of length (at least) $c(\log(\ep^{-1}))^{-2}$ which will be very close to 
a straight line, see figure \ref{pic}. Let $s<0$, $r>0$ be such that $X(s),X(e)$ are the endpoints of the path 
(where we assume without loss of generality $X(s)$ is closer to $\partial B_1(0)$ than $X(e)$). If we are able to 
show that $X(s)\in \partial T_1\cap \partial B_1(0)$ then the argument can proceed very much as described in the paragraphs above. The only way this can fail is if the path is (close to) a line of length $c(\log(\ep^{-1}))^{-1}$ and 
runs, (roughly speaking) parallel to $\partial T_1\cap \partial B_1(0)$. However as 
$\lt|u(X(e))-u(X(s))\rt|\geq c(\log(\ep^{-1}))^{-1}$ this implies we must have 
$\lt|u(X(e))\rt|\geq c(\log(\ep^{-1}))^{-1}$, but since the path is close to `parallel' to $\partial B_1(0)\cap 
\partial T_1$ we have $\mathrm{dist}(X(e),\partial B_1(0))\leq c\log(\ep^{-1})^{-2}$ 
which contradicts $1$-Lipschitz type property as represented by inequality (\ref{avvu5}), thus we must have 
that $X(s)\in \partial T_1\cap \partial B_1(0)$. By use of this argument assumption (\ref{fa1}) can be avoided.

\section{The E.L. equation}

Note that if $u$ is a critical point of $I_{\ep}$ it weakly satisfies the E.L. equation i.e.\
\begin{equation}
\label{reb6}
\ep\Delta^2u+\ep^{-1}\divv\lt(\lt(1-\lt|Du\rt|^2\rt)Du\rt)=0.
\end{equation}
Let $w\in W^{1,1}$ define $w_i:=\frac{\partial w}{\partial x_i}$, similarly for 
$v\in W^{2,1}, s\in W^{3,1}$ define $v_{ij}:=\frac{\partial^2 v}{\partial x_i \partial x_j}$ and 
$s_{ijk}:=\frac{\partial^3 s}{\partial x_i\partial x_j\partial x_k}$.

\begin{a1}
\label{RL1} Suppose $u\in W^{2,2}(\Omega)$ is a weak solution of (\ref{reb6}). Define $\Omega_{\ep^{-1}}:=\ep^{-1}\Omega$ and let
$v:\Omega_{\ep^{-1}}\rightarrow \R$ be defined by
$v\lt(z\rt):=u\lt(\ep z\rt)\ep^{-1}$, then $v$ satisfies
\begin{equation}
\label{le151}
\Delta^2 v+\divv\lt(\lt(1-\lt|Dv\rt|^2\rt)Dv\rt)=0
\end{equation}
weakly in $\Omega_{\ep^{-1}}$.
\end{a1}
\em Proof. \rm Follows directly from the definition of $u$. 

%
%

%
%
%
%
%
%

\begin{a1}
\label{RL2} We will show that any $v\in
W^{2,2}\lt(\Omega_{\ep^{-1}}\rt)$ that satisfies (\ref{le151})
weakly in $\Omega_{\ep^{-1}}$ is such that for any $U\subset\subset \Omega_{\ep^{-1}}$, $v\in
W^{3,2}\lt(U\rt)$ and $v$ satisfies
\begin{equation}
\label{reb1} \int \sum_{i,j,p=1}^2 v_{ijp}\phi_{ijp}+\lt(\lt(1-\lt|Dv\rt|^2\rt)\cdot Dv\rt)_p D\phi_p \;dz=0
\end{equation}
for any $\phi\in C^1_0\lt(U\rt)$.

\end{a1}
\em Proof. \rm Given set $S\subset \R^2$, let 
$d(x,S)=\inf\lt\{\lt|z-x\rt|:z\in S\rt\}$ and define $N_{\delta}(S):=\lt\{x:d(x,S)<\delta\rt\}$.

Step 1. For $\delta>0$ let $\Pi_{\delta}:=\Omega_{\ep^{-1}}\backslash N_{\delta}(\partial \Omega_{\ep^{-1}})$. 
We will show that $D^2 v\in W^{1,2}(\Pi_{3\delta})$. 

\em Proof of Step 1. \rm Let $g(x):=Dv(x)\lt(1-\lt|Dv(x)\rt|^2\rt)$ and $w:=\Delta v$. Since 
$v\in W^{2,2}(\Omega_{\ep^{-1}})$, by Poincare's inequality (Theorem 2, Section 4.5.2 \cite{evans3}) $Dv\in L^p(\Omega_{\ep^{-1}})$ for any $p<\infty$, hence $g\in L^q(\Omega_{\ep^{-1}})$ for any $q<\infty$. So
$$
\int w\Delta \phi=\int g \cdot D\phi\text{ for any }\phi\in C_0^{\infty}(\Omega_{\ep^{-1}}). 
$$ 

Let $\rho\in C_0^{\infty}(B_1)$ be the standard convolution kernel and define 
$\rho_{\sigma}(z)=\rho\lt(\frac{z}{\sigma}\rt)\sigma^{-2}$. Given function $f\in W^{1,1}$ we 
denote the convolution of $f$ and $\rho_{\sigma}$ by $f*\rho_{\sigma}$. Let $\veps\in (0,\delta)$ and define $w_{\veps}:=w*\rho_{\veps}$ and $g_{\veps}:=g*\rho_{\veps}$. Now for any $\phi\in C_0^{\infty}(\Omega_{\ep^{-1}})$, 
defining $\phi_{\veps}=\phi*\rho_{\veps}$  we have 
$$
\int w_{\veps}\Delta \phi=\int w\Delta \phi_{\veps}=\int g \cdot D\phi_{\veps}=\int g_{\veps} \cdot D\phi
$$
which gives that $\Delta w_{\veps}(z)=-\mathrm{div} g_{\veps}(z)$ for any $z\in \Pi_{\delta}$. 
Let $\psi\in C_0^{\infty}(\Pi_{\delta})$ with $\psi=1$ on $\Pi_{2\delta}$ and $\lt|D \psi\rt|<c\delta^{-1}$ and 
$\lt|D^2\psi\rt|<c\delta^{-2}$. Define $s(x)=w_{\veps}(x)\psi(x)$, so 
$$
\Delta s=-\mathrm{div} g_{\veps}\psi+2 Dw_{\veps}\cdot D\psi+w_{\veps}\Delta \psi.
$$ 
Now $\mathrm{div}(g_{\veps}\psi)=\mathrm{div}g_{\veps}\psi+g_{\veps}\cdot D\psi$ and 
$2 Dw_{\veps}\cdot D\psi=\mathrm{div}(2 w_{\veps} D\psi)-2w_{\veps}\Delta \psi$ and thus 
\begin{equation}
\label{xc1}
\Delta s=\mathrm{div}(-g_{\veps}\psi+2 w_{\veps}D\psi)+g_{\veps}\cdot D\psi- w_{\veps} \Delta\psi.
\end{equation} 
Let $X=Ds$, so by (\ref{xc1}) we have that 
\begin{equation}
\label{xc2}
\mathrm{curl}(X)=0\text{ and }\mathrm{div}(X+g_{\veps}\psi-2 w_{\veps}D\psi)=g_{\veps}\cdot D\psi-w_{\veps}\Delta \psi.
\end{equation} 

For any $C^2$ vector field $V$, let $H(V)$ denote the Hodge projection of $V$ onto the subspace of curl free vector fields, 
i.e.\ $H(V)=-D \Delta^{-1} \mathrm{div} V$, so $H(V)$ satisfies $\mathrm{div}(H(V)+V)=0$ and $\mathrm{curl}H(V)=0$ on 
$\R^2$. So from (\ref{xc2}) then we have 
\begin{equation}
\label{xc2.5}
\mathrm{curl}(X-H(g_{\veps}\psi-2 w_{\veps}D\psi))=0\text{ and }\mathrm{div}(X-H(g_{\veps}\psi-2 w_{\veps}D\psi))
=g_{\veps}\cdot D\psi-w_{\veps}\Delta \psi.
\end{equation}
Let $\eta\in C^{\infty}(\R^2)$ be such that 
\begin{equation}
\label{ulk3.3}
D\eta=X-H(g_{\veps}\psi-2 w_{\veps}D\psi), 
\end{equation}
so finally we have 
\begin{equation}
\label{ulk2}
\Delta \eta=g_{\veps}\cdot D\psi-w_{\veps}\Delta\psi.
\end{equation}
Now recall $X=Ds$ where $s=w_{\veps}\psi$. Thus $D s=D w_{\veps}\psi+w_{\veps}D\psi$ and thus for any $p\in \lt[1,2\rt]$,
\begin{equation}
\label{ggfeq2}
\|X\|_{L^p(\R^2)}\leq c\|D w_{\veps}\|_{L^p(\R^2)}+c\|w_{\veps}\|_{L^p(\R^2)}\leq c\|w*D\rho_{\veps}\|_{L^p(\R^2)}+c\|w_{\veps}\|_{L^p(\R^2)}\leq 
c\veps^{\frac{2-3p}{p}}\|D^2 u\|_{L^2(\Omega_{\ep^{-1}})}\leq c\veps^{\frac{2-3p}{p}}.
\end{equation}
And by $L^p$ boundedness of Hodge projection we know 
\begin{equation}
\label{ggfeq1}
\|H(g_{\veps}\psi-2w_{\veps}D\psi)\|_{L^p(\R^2)}\leq c\|g_{\veps}\psi-2w_{\veps}D\psi\|_{L^p(\R^2)}\leq 
c\|g_{\veps}\|_{L^p(\Omega_{\ep^{-1}})}+c\|w_{\veps}\|_{L^p(\Omega_{\ep^{-1}})}\leq c.
\end{equation}
Thus for $p=\frac{3}{2}$ we have $\|D\eta\|_{L^{\frac{3}{2}}(\R^2)}\overset{(\ref{ggfeq1}),(\ref{ggfeq2}),(\ref{ulk3.3})}{\leq} 
c\veps^{-\frac{5}{3}}$. 
What we need to do is obtain an $\veps$ independent bound on $D\eta$, we will achieve this by use of (\ref{ulk2}). First note by Holder 
$g_{\veps}\cdot D\psi-w_{\veps}\Delta \psi\in L^{\frac{3}{2}}(\R^2)$ from (\ref{ulk2}) by Standard $L^p$ estimates on Riesz transforms (see 
Proposition 3, Section 1.3. Chapter 3 \cite{stein}) we know 
\begin{equation}
\label{ggfeq1.6}
\|D^2 \eta \|_{L^{\frac{3}{2}}(\R^2)}\leq c \|g_{\veps} \|_{L^{\frac{3}{2}}(\Omega_{\ep^{-1}})}
+c\|w_{\veps}\|_{L^{\frac{3}{2}}(\Omega_{\ep^{-1}})}\leq c.
\end{equation}
So $D\eta \in W^{1,\frac{3}{2}}(\R^2)$ and thus by Sobolev embedding theorem (Theorem 1, Section 4.5.1. \cite{evans3}) we have 
$\|D\eta\|_{L^6(\R^2)}\leq c\|D^2 \eta\|_{L^{\frac{3}{2}}(\R^2)}\overset{(\ref{ggfeq1.6})}{\leq} c$. 
As $\spt X\subset \Pi_{\delta}\subset \Omega_{\ep^{-1}}$, $\|D s\|_{L^2(\R^2)}=\|D s\|_{L^2(\Omega_{\ep^{-1}})}\leq c$ and using $L^2$ boundedness of the Hodge projection 
\begin{equation}
\label{ulk3}
\|Ds\|_{L^2(\R^2)}\overset{(\ref{ulk3.3})}{\leq} \|D\eta\|_{L^2(\Omega_{\ep^{-1}})}
+\|H(g_{\veps}\psi-2 w_{\veps} D\psi)\|_{L^2(\Omega_{\ep^{-1}})}\leq c.
\end{equation}
Since $Ds=Dw_{\veps}\psi+w_{\veps} D\psi$, so $\|Dw_{\veps}\psi\|_{L^2(\R^2)}
\overset{(\ref{ulk3})}{\leq} c
+\|w_{\veps}D\psi\|_{L^2(\R^2)}$. Now $w_{\veps}=\triangle v_{\veps}$ and so 
$\|w_{\veps}D\psi\|_{L^2(\R^2)}\leq c\|D^2 v_{\veps}\|_{L^2(\Pi_{\delta})}\leq c$ for any $\veps>0$. 
Hence 
\begin{equation}
\label{finqeq1}
\|D w_{\veps}\|_{L^2(\Pi_{2\delta})}<c \text{ for all }\veps>0.
\end{equation}

Let $q\in C^{\infty}_0\lt(\Pi_{2\delta}\rt)$ with 
$q\equiv 1$ on $\Pi_{3\delta}$. Let 
$z_{\veps}=v_{\veps,1}q$ so $\triangle z_{\veps}=\triangle v_{\veps,1} q+2Dv_{\veps,1}\cdot Dq+v_{\veps,1}\triangle q$. Thus 
as $\triangle v_{\veps,1}=w_{\veps,1}$
$$
\|\triangle z_{\veps}\|_{L^2(\R^2)}\leq 
\|\triangle v_{\veps,1}q\|_{L^2(\R^2)}
+2\|D v_{\veps,1}\cdot Dq\|_{L^2(\R^2)}+
\|v_{\veps,1}\triangle q\|_{L^2(\R^2)} \overset{(\ref{finqeq1})}{\leq} c. 
$$
Now as we have seen before by $L^2$ estimates on Riesz transforms, this implies $D^2 z_{\veps}\in L^2(\R^2)$. As $D^2 z_{\veps}=D^2 v_{\veps,1} q+2 Dv_{\veps,1}\otimes Dq+v_{\veps,1}D^2 q$ we have that 
\begin{equation}
\int_{\Pi_{3\delta}} \lt|D^2 v_{\veps,1}\rt|^2 dx\leq c\int_{\R^2} \lt|D^2 z_{\veps}\rt|^2 dx+c\int_{\R^2} \lt|D v_{\veps,1}\rt|^2+c\int_{\R^2} \lt|v_{\veps,1}\rt|^2 dx\leq c\text{ for every }\veps>0.
\end{equation}
Arguing in exactly the same way gives $\int_{\Pi_{3\delta}} \lt|D^2 v_{\veps,2}\rt|^2 dx\leq c$ for every $\veps>0$, thus 
$$
\int_{\Pi_{3\delta}} \lt|D^3 v_{\veps}\rt|^2\leq c\text{ for every }\veps>0.
$$

Now for any $\veps_n\rightarrow 0$, $D^2 v_{\veps_n}$ is a bounded sequence in 
$W^{1,2}(\Pi_{3\delta})$, so for some subsequence 
$k_n$, $D^2 v_{\veps_{k_n}}\rightharpoonup \zeta\in W^{1,2}(\Pi_{3\delta}:\R^{2\times 2})$. Clearly $\zeta=D^2 v$ for a.e.\ in $\Pi_{3\delta}$.  
Let $i,j,k\in\lt\{1,2\rt\}$ and $\phi\in C_0^{\infty}\lt(\Pi_3\rt)$, 
\begin{eqnarray}
\int v_{,ij} \phi_{,k}&=&\lim_{n\rightarrow \infty} \int v_{\veps_{k_n},ij} \phi_{,k} dx\nn\\
&=&\lim_{n\rightarrow \infty} \int -v_{\veps_{k_n},ijk} \phi dx\nn\\
&=&\int -\zeta_{ij,k} \phi dx.\nn
\end{eqnarray}
Thus $v_{,ij}\in W^{1,2}(\Pi_{3\delta})$ for any $i,j\in \lt\{1,2\rt\}$ and hence 
$D^2 v\in W^{1,2}(\Pi_{3\delta})$.\nl

\em Step 2. \rm We will show that $v$ satisfies (\ref{reb1}). 

\em Proof of Step 2. \rm Take any arbitrary $\phi\in C^{\infty}(\Omega_{\ep^{-1}})$, letting $\psi^{h}(z):=\frac{\phi(z+he_p)-\phi(z)}{h}$ we know from (\ref{le151})
\begin{eqnarray}
\label{reb41}
&~&
\int \sum_{i,j} v_{ij}\lt(y\rt)\phi_{ijp}\lt(y\rt)+
\lt(1-\lt|Dv\lt(y\rt)\rt|^2\rt)Dv\lt(y\rt) D\phi_p\lt(y\rt) dy\nn\\
&~&\qd\qd\qd=\lim_{h\rightarrow 0}h^{-1}\int \sum_{i,j=1}^{2}
v_{ij}\lt(y\rt)\psi^h_{ij}\lt(y\rt)
+\lt(1+\lt|Dv\lt(y\rt)\rt|^2\rt)Dv\lt(y\rt)D\psi^h\lt(y\rt) dy\nn\\
&~&\qd\qd\qd =0
\end{eqnarray}
thus integrating by parts
$$
\int \sum_{i,j} v_{ijp}\phi_{ij}+\lt(\lt(1-\lt|Dv\rt|^2\rt)Dv\rt)_p D\phi dy=0.
$$
Repeating the argument gives us (\ref{reb1}). \;\;\;\;$\Box$

%
%
%

\begin{a1}
\label{L0.5} Let $u\in W_0^{2,2}(B_1(0))$ be the minimizer of $I_{\ep}$, then 
\begin{equation}
\label{uuzu5}
I_{\ep}(u)\leq c\ep\log(\ep^{-1}).
\end{equation}
\end{a1}
\em Proof. \rm Let $\rho$ be the standard rotationally symmetric convolution kernel 
with $\spt \rho\subset B_{2}(0)$ and let $\rho_{\ep}(z):=\rho(\frac{z}{\ep})\ep^{-2}$. Let 
$w(x)=1-\lt|x\rt|$ and $w_{\ep}=w*\rho_{\ep}$. So if $y\in B_{4\ep}(0)$
\begin{equation}
\label{fa2}
\lt| D^2 w_{\ep}(y)\rt|\leq \lt|\int (w(z)-1)D^2 \rho_{\ep}(y-z) dz\rt|\leq 
c\ep^{-4}\int_{B_{6\ep}(0)} \lt|w(z)-1\rt| dz\leq c\ep^{-1}.
\end{equation}

Note $D w(y)=-\frac{y}{\lt|y\rt|}$ and $D^2 w(y)=\frac{y\otimes y}{\lt|y\rt|^3}-\lt|y\rt|^{-1}Id$ so 
$\lt|D^2 w(y)\rt|\leq \frac{4}{\lt|y\rt|}$. So 
\begin{equation}
\label{fa3}
\lt|D^2 w_{\ep}(y)\rt|\leq \lt|\int D^2 w(z)\rho_{\ep}(y-z) dz\rt|\leq 4\int \frac{\rho_{\ep}(y-z)}{\lt|z\rt|} dz\leq 
\frac{c}{\lt|y\rt|} \text{ for any } y\not\in B_{4\ep}(0) .
\end{equation}
Thus 
$$
\int_{B_1(0)} \lt|D^2 w_{\ep}\rt|^2 dy\leq \int_{B_{4\ep}(0)} \lt|D^2 w_{\ep}\rt|^2 dy
+\int_{B_1(0)\backslash B_{4\ep}(0)} \lt|D^2 w_{\ep}\rt|^2 dy\overset{(\ref{fa2}),(\ref{fa3})}{\leq} c+c\int_{4\ep}^1 r^{-1} dr\leq c\log(\ep^{-1}).
$$

Now $\lt\{x\in \R^2:w_{\ep}(x)=0\rt\}$ is a circle of radius $h\eqsim 1$ so defining 
$v(x)=w_{\ep}\lt(\frac{x}{h}\rt)h$, $v\in W_0^{2,2}(B_1(0))$ and 
$\int_{B_1(0)} \lt|D^2 v\rt|^2 dx\leq c\log(\ep^{-1})$. Now if $x\not \in B_{4\ep}(0)$, 
$\lt|D w_{\ep}(x)-D w(x)\rt|=\lt|\int (D w(z)-D w(x))\rho_{\ep}(x-z) dz\rt|\leq \frac{c\ep}{\lt|x\rt|}$. So 
$\lt|\lt|D w_{\ep}(x)\rt|^2-1\rt|^2\leq c\lt|\lt|D w_{\ep}(x)\rt|-1\rt|^2\leq \frac{c\ep^2}{\lt|x\rt|^2}$. Thus 
\begin{eqnarray}
\int_{B_1(0)} \lt|1-\lt|D w_{\ep}(x)\rt|^2\rt|^2 dx&\leq& c\ep^2+\int_{B_1(0)\backslash B_{4\ep}(0)} 
\lt|1-\lt|D w_{\ep}(x)\rt|^2\rt|^2 dx\nn\\
&\leq&c\ep^2+\int_{4\ep}^1 \frac{\ep^2}{r} dr\nn\\
&\leq& c\log(\ep^{-1})\ep^2\nn
\end{eqnarray}
 and this establishes (\ref{uuzu5}). $\Box$

%
%
%

\begin{a1}
\label{L2} Let $u\in W^{2,2}_0(B_1(0))$ be a minimizer of $I_{\ep}$. Let $\cv_1$ be a some small positive constant to be 
chosen later. Define $A(x,\alpha,\beta):=B_{\beta}(x)\backslash \overline{B_{\alpha}(x)}$. We divide
$B_{1}(0)$ into $N=\lt[\cv_1^{-2}\log(\ep^{-1})\rt]$ slices of equal angle, denote their closure by
$T_1,T_2,\dots T_N$. There must exists a set $\Pi\subset\lt\{1,2,\dots N\rt\}$ with 
$\ca{\Pi}\geq \frac{N}{2}$ such that if $i\in \Pi$
\begin{eqnarray}
\label{avv14.7}
&~&\inf\lt\{\lt|Du\lt(z\rt)\rt|:z\in  T_i\cap A(0,c\log(\ep^{-1})\ep,1-2\ep)\rt\}
>\frac{1}{2}\text{ and }\nn\\
&~&\qd\qd
\sup\lt\{\lt|Du\lt(z\rt)\rt|:z\in T_i\cap A(0,c\log(\ep^{-1})\ep,1-2\ep)   \rt\}<2.
\end{eqnarray}
\end{a1}

\em Proof of Lemma \ref{L2}. \rm Define $v\lt(z\rt)=u\lt(\ep z\rt)\ep^{-1}$. Let $S_i=\ep^{-1} T_i$ for $i=1,2, \dots N$. 
For $i\in\lt\{2,3,\dots N-1\rt\}$ define
$$
\wt{S}_i=S_{i-1}\cup S_i\cup S_{i+1}\text{ and let }
\wt{S}_1=S_{N-1}\cup S_1\cup S_{2},\; \wt{S}_N=S_{N-1}\cup S_N\cup S_{1}.
$$
Define 
\begin{equation}
\label{av1} G_0:=\lt\{i\in\lt\{1,2,\dots N\rt\}:\int_{\wt{S}_i}
\lt|1-\lt|Dv\rt|^2\rt|^2+\lt|D^2 v\rt|^2 dz\leq \cv_1\rt\}.
\end{equation}
Note that by (\ref{uuzu5}) of Lemma \ref{L0.5} we know 
$\int_{B_{\ep^{-1}}(0)} \lt|1-\lt|Dv\rt|^2\rt|^2+\lt|D^2 v\rt|^2 dx\leq c\log(\ep^{-1})$, so $\cv_1(N-\ca{G_0})\leq c\log(\ep^{-1})$, thus 
(assuming we chose $\cv_1$ small enough) $\frac{\cv_1^{-2}}{2}\log(\ep^{-1})\leq \ca{G_0}$.\nl

\em Step 1. \rm Let $i\in G_0$, we will show that for any $y_0\in \wt{S}_i$ such that 
$B_2\lt(y_0\rt)\subset \wt{S}_i$ and $\psi\in C^{\infty}_0\lt(B_2\lt(y_0\rt)\rt)$ such that $\psi\equiv
1$ on $B_1\lt(y_0\rt)$ we have 

\begin{equation}
\label{reb11} \int \lt|D^3 v\rt|^2\psi^6 dz\leq c.
\end{equation}

\em Proof of Step 1. \rm Let
$Y=\lt(4\pi\rt)^{-1}\int_{B_2\lt(y_0\rt)} Dv$,
$T=\lt(4\pi\rt)^{-1}\int_{B_2\lt(y_0\rt)} v$ and
we define $\ti{v}\lt(z\rt)=v\lt(z\rt)-Y\cdot\lt(z-y_0\rt)-T$.

Let $\phi:=\ti{v}\psi^6$. So
$\phi_p=\ti{v}_p\psi^6+6\ti{v}\psi^5\psi_p$ and
\begin{equation}
\label{reb8.8}
\phi_{pi}=v_{pi}\psi^6+6\ti{v}_p\psi^5\psi_i+6\ti{v}_i\psi^5\psi_p+6\ti{v}\lt(\psi^5\psi_p\rt)_i.
\end{equation}
\begin{eqnarray}
\label{reb9}
\phi_{pij}&=&v_{pij}\psi^6+6 v_{pi}\psi^5\psi_j+6 v_{pj}\psi^5\psi_i+6\ti{v}_p\lt(\psi^5\psi_i\rt)_j\nn\\
&~&+6v_{ij}\psi^5\psi_p+6\ti{v}_i\lt(\psi^5\psi_p\rt)_j
+6\ti{v}_j\lt(\psi^5\psi_p\rt)_i+6\ti{v}\lt(\psi^5\psi_p\rt)_{ij}.
\end{eqnarray}

By the fact that $B_2(y_0)\subset \wt{S}_i$ we know $\int_{B_2(y_0)} \lt|D^2 v\rt|^2 \leq \cv_1$, by 
Poincare's inequality this implies $\|D\ti{v}\|_{L^2(B_2(y_0))}\leq c$ and $\|\ti{v}\|_{L^2(B_2(y_0))}\leq c$. So 
from (\ref{reb9})
\begin{eqnarray}
\label{reb37} \lt|\int v_{ijp}\phi_{ijp}-\int
\lt(v_{ijp}\rt)^2\psi^6 \rt|&\overset{(\ref{reb9})}{\leq}&c\|v_{ijp}\psi^3\|_{L^2}\lt(\|D^2
v\|_{L^2\lt(B_2(y_0)\rt)}+\|D
\ti{v}\|_{L^2\lt(B_2(y_0)\rt)}+
\|\ti{v}\|_{L^2\lt(B_2(y_0)\rt)}\rt)\nn\\
&\leq &c\|D^3 v \psi^3\|_{L^2}.
\end{eqnarray}
Now 
\begin{eqnarray}
\label{reb23} \lt|\int \lt(\lt(1-\lt|D v\rt|^2\rt)D v\rt)_p\cdot
D\phi_p\; dz\rt|&=&\lt|\int \lt(\lt(1-\lt|D v\rt|^2\rt)D v\rt)\cdot
D\phi_{pp} dz\rt|\nn\\
&\leq&\lt|\int \lt(\lt(1-\lt|D v\rt|^2\rt)D v\rt)\cdot
\lt(D\phi_{pp}-Dv_{pp}\psi^6\rt) dz\rt|\nn\\
&~&+\lt|\int \lt(\lt(1-\lt|D v\rt|^2\rt)D v\rt)\cdot D
v_{pp}\psi^6 dz\rt|\nn\\
&\overset{(\ref{reb9})}{\leq}&
c\|\lt(1-\lt|Dv\rt|^2\rt)Dv\|_{L^2\lt(B_2\lt(y_0\rt)\rt)}\|D^2
v\|_{L^2\lt(B_2\lt(y_0\rt)\rt)}\nn\\
&~&+\|D^3 v\psi^3\|_{L^2}\|\lt(1-\lt|Dv\rt|^2\rt)Dv\psi^3\|_{L^2}\nn\\
&\overset{(\ref{av1})}{\leq}&c\lt(1+\|D^3 v \psi^3\|_{L^2\lt(B_2\lt(y_0\rt)\rt)}\rt).
\end{eqnarray}

Recalling the fact that by Lemma \ref{RL2}, $v$ satisfies (\ref{reb1}) we have
\begin{eqnarray}
\lt|\int \sum_{i,j,p=1}^{2}
\lt(v_{ijp}\rt)^2\psi^6 dz\rt|&\overset{(\ref{reb1})}{=}&
\lt|\int\sum_{i,j,p=1}^{2}
\lt(v_{ijp}\rt)^2\psi^6-v_{ijp}\phi_{ijp}-\int\lt(\lt(1-\lt|Dv\rt|^2\rt)Dv\rt)_p\cdot
D\phi_p dz\rt|\nn\\
&\overset{(\ref{reb37}),(\ref{reb23})}{\leq}&c\|D^3v\psi^3\|_{L^2}+c.\nn
\end{eqnarray}
And this establishes (\ref{reb11}).\nl

\em Proof of Lemma \ref{L2}. \rm  By Theorem 2, Section 5.6 \cite{evans2}
$$
\|D^2 v\|_{L^4(B_2(y_0))}\leq \|D^2 v\|_{W^{1,2}(B_2(y_0))}\leq c+\|D^3 v\|_{L^2(B_2(y_0))}
\overset{(\ref{reb1})}{\leq} c.
$$
By Sobolev embedding this implies $Dv$ is $\frac{1}{2}$-Holder in $B_1\lt(y_0\rt)$.

Since $\int_{B_1\lt(y_0\rt)} \lt|1-\lt|Dv\rt|^2\rt|^2 dz\leq \cv_1$. Let 
$L=\lt\{z\in B_1\lt(y_0\rt):\lt|1-\lt|Dv\rt|^2\rt|^2\leq
\sqrt{\cv_1}\rt\}$
so we have $\lt|B_1\lt(y_0\rt)\backslash L\rt|\leq
\sqrt{\cv_1}$. So $B_{4\cv_1^{\frac{1}{4}}}\lt(y_0\rt)\cap L\not=\emptyset$ so
we can pick $z_1\in B_{4 \cv_1^{\frac{1}{4}}}(y_0)\cap L$. Since $Dv$ is $\frac{1}{2}$ Holder
\begin{eqnarray}
\lt|\lt|Dv\lt(y_0\rt)\rt|-1\rt|&\leq&\lt|Dv\lt(y_0\rt)-Dv\lt(z_1\rt)\rt|+\cv_1^{\frac{1}{4}}\nn\\
&\leq&c\lt|y_0-z_1\rt|^{\frac{1}{2}}+\cv_1^{\frac{1}{4}}\nn\\
&\leq&c\cv_1^{\frac{1}{8}},\nn
\end{eqnarray}
assuming we chose $\cv_1$ small enough this implies $\lt|Dv(y_0)\rt|\in (\frac{1}{2},2)$. Since 
$y_0$ is an arbitrary point in $\wt{S_i}\backslash N_2(\partial \wt{S_i})$ and $Du(\ep y_0)=Dv(y_0)$ this 
implies (\ref{avv14.7}).  $\Box$

%
%
%
%

\begin{a1}
\label{L5} Let $u\in W^{2,2}(B_1(0))$. 
Suppose 
\begin{equation}
\label{qza1}
\int_{B_1(0)} \lt|1-\lt|Du\rt|^2\rt|\lt|D^2 u\rt| dz\leq \beta
\end{equation}
and 
\begin{equation}
\label{qza2}
\int_{B_1(0)} \lt|1-\lt|Du\rt|^2\rt| dz\leq \beta.
\end{equation}

We will show that for any $w\in S^1$ we can find a set $G_w\subset P_{w^{\perp}}\lt(B_1(0)\rt)$ with
\begin{equation}
\label{avv14.5}
\lt|P_{w^{\perp}}\lt(B_1(0)\rt)\backslash G_{w}\rt|\leq
\beta^{\frac{1}{3}}
\end{equation}
and for any $x\in G_w$ we have
\begin{equation}
\label{avv14} \sup\lt\{\lt|\lt|Du\lt(z\rt)\rt|-1\rt|:z\in
P^{-1}_{w^{\perp}}\lt(x\rt)\cap B_1(0)\rt\}\leq
5\beta^{\frac{1}{3}}.
\end{equation}
\end{a1}

\em Proof of Lemma. \rm Let
$$
B_w:=\lt\{x\in P_{w^{\perp}}\lt(B_1(0)\rt):
\int_{P^{-1}_{w^{\perp}}\lt(x\rt)\cap B_1(0)}
\lt|1-\lt|Du\rt|^2\rt|\lt|D^2 u\rt|+\lt|1-\lt|Du\rt|^2\rt| dz\leq
\beta^{\frac{2}{3}}\rt\}.
$$
By Chebyshev's inequality we have
$\lt|P_{w^{\perp}}\lt(B_1(0)\rt)\backslash
B_{w}\rt|\leq 2\beta^{\frac{1}{3}}$. For any $x\in P_{w^{\perp}}(B_{1-\beta^{\frac{2}{3}}}(0))$ we know 
$\lt|P^{-1}_{w^{\perp}}(x)\cap B_1(0)\rt|\geq \beta^{\frac{1}{3}}$ and so if in addition $x\in B_w$ we have 
that there must exists $z_x\in P^{-1}_{w^{\perp}}(x)\cap B_1(0)$ such that 
$\lt|1-\lt|Du(z_x)\rt|\rt|\leq \beta^{\frac{1}{3}}$. 

Suppose $x\in B_w\cap P_{w^{\perp}}(B_{1-\beta^{\frac{2}{3}}}(0))$ and 
for some $y_x\in P^{-1}_{w^{\perp}}(x)\cap B_1(0)$ we have $\lt|1-\lt|Du(y_x)\rt|\rt|\geq 5\beta^{\frac{1}{3}}$. 
Then as we can assume without loss of generality that $Du$ is continuous on $P^{-1}_{w^{\perp}}(x)\cap B_1(0)$ and 
so there must exists $a_x, b_x\in P^{-1}_{w^{\perp}}(x)\cap B_1(0)$ such that $\lt|\lt|Du(a_x)\rt|-\lt|Du(b_x)\rt|\rt|\geq 
\beta^{\frac{1}{3}}$ and $\inf\lt\{\lt|Du(x)\rt|:x\in [a_x,b_x]\rt\}\geq 1+4\beta^{\frac{1}{3}}$. However by the 
fundamental theorem of Calculus 
$$
4\beta^{\frac{1}{3}}\lt|\lt|Du(a_x)\rt|-\lt|Du(b_x)\rt|\rt| \leq \int_{a_x}^{b_x} \lt|1-\lt|Du\rt|\rt|\lt|D^2 u\rt|
\leq\beta^{\frac{2}{3}}
$$
which is a contradiction. Thus taking $G_w:=B_w\cap P_{w^{\perp}}(B_{1-\beta^{\frac{1}{3}}}(0))$ completes the proof 
of the lemma. $\Box$\nl\nl

%
%
%

\begin{a1}
\label{LL5}

Suppose $\ui$ is a $C^2$ function that satisfies (\ref{qza1}), (\ref{qza2}) and $\Lambda\subset B_1(0)$ is convex with the property that 
$\inf\lt\{\lt|D\ui(x)\rt|:x\in \Lambda\rt\}>\frac{1}{3}$ and  $\sup\lt\{\lt|D\ui(x)\rt|:x\in \Lambda\rt\}<3$. 

Given function $X:\R\rightarrow \R^2$ that solves $X(0)=x$ and $\dot{X}(s)=D\ui(X(s))$, 
suppose $s_1<0<s_2$ are such that $X(s)\in \Lambda$ for any $s\in [s_1,s_2]$  then 
\begin{equation}
\label{zaq4}
\ui(X(s_2))-\ui(X(s_1))\geq (1-\beta^{\frac{1}{3}})\lt|X(s_2)-X(s_1)\rt|-c\beta^{\frac{1}{3}}.
\end{equation}
And if in addition $X(s_1),X(s_2)\not\in B_r(x)$ for some $B_r(x)\subset \Omega$, then  
\begin{equation}
\label{zaq4.5}
\lt\{X(s):s\in [s_1,s_2]\rt\}\subset N_{c\frac{\beta^{\frac{1}{6}}}{\sqrt{r}}}(\lt[X(s_1),X(s_2)\rt]).
\end{equation}
\end{a1}

\em Proof. \rm Let $w\in S^1$ be orthogonal to $X(s_2)-X(s_1)$. Let $G_w$ be the set satisfying (\ref{avv14.5}) and (\ref{avv14}) 
from Lemma \ref{L5}. Let $P=\lt\{X(t):t\in\lt[s_1,s_2\rt]\rt\}$ and 
$\Gamma=P\cap P_{w^{\perp}}^{-1}(G_w)$. So $H^1(\Gamma)\geq \lt|P_{w^{\perp}}(\lt[X(s_1),X(s_1)\rt])\cap G_w\rt|\geq \lt|X(s_2)-X(s_1)\rt|-\beta^{\frac{1}{3}}$ and so 
\begin{eqnarray}
\label{vaq8}
\ui(X(s_2))-\ui(X(s_1))&=&\int_{P} D\ui(z)\cdot t_z dH^1 z\nn\\
&\geq& (1-c\beta^{\frac{1}{3}})H^1(\Gamma)+\frac{1}{3}H^1(P\backslash \Gamma)\nn\\
&\geq&(1-c\beta^{\frac{1}{3}})\lt|X(s_2)-X(s_1)\rt|+\frac{1}{3}H^1(P\backslash \Gamma)-c\beta^{\frac{1}{3}}
\end{eqnarray}
which establishes (\ref{zaq4}). Now 
\begin{eqnarray}
\label{vaq9}
\ui(X(s_2))-\ui(X(s_1))&\leq& \int_{\lt[X(s_1),X(s_2)\rt]} \lt|D\ui(z)\rt| dH^1 z\nn\\
&\leq& (1+c\beta^{\frac{1}{3}})\lt|P_{v^{\perp}}(\lt[X(s_2),X(s_1)\rt]\cap G_w)\rt|+
3\lt|P_{v^{\perp}}(\lt[X(s_2),X(s_1)\rt]\backslash G_w)\rt|\nn\\
&\leq &\lt|X(s_2)-X(s_1)\rt|+c\beta^{\frac{1}{3}}
\end{eqnarray}
now putting (\ref{vaq8}) and (\ref{vaq9}) together we have $H^1(P\backslash \Gamma)\leq c\beta^{\frac{1}{3}}$. Now this 
and the second inequality of (\ref{vaq8}) and inequality (\ref{vaq9}) imply that 
\begin{equation}
\label{ulm1}
\lt|X(s_2)-X(s_1)\rt|-c\beta^{\frac{1}{3}}\geq H^1(P). 
\end{equation}
If 
$X(s_1),X(s_2)\not\in B_r(x)$ then as $X(0)=x\in P$ and as $P$ is connected we know 
$H^1(P)\geq \lt|X(s_1)-X(0)\rt|+\lt|X(s_2)-X(0)\rt|\geq 2 r$ which by (\ref{ulm1}) implies  
$\lt|X(s_1)-X(s_2)\rt|\geq r$ and so $\lt|X(s_1)-X(s_2)\rt|(1+\frac{c\beta^{\frac{1}{3}}}{r})\geq H^1(P)$.  
Now letting $t_z$ denote the tangent to the curve $P$ at point $z$ we have 
\begin{eqnarray}
\int_{P} \lt|t_z-\frac{X(s_2)-X(s_1)}{\lt|X(s_2)-X(s_1)\rt|}\rt|^2 dH^1 z&=&
\int_{P} 2-2t_z\cdot\lt(\frac{X(s_2)-X(s_1)}{\lt|X(s_2)-X(s_1)\rt|}\rt) dH^1 z\nn\\
&=&2H^1(P)-2\lt|X(s_2)-X(s_1)\rt|\nn\\
&\leq&\frac{c\beta^{\frac{1}{3}}}{r}.\nn
\end{eqnarray}
By Holder's inequality and the fundamental theorem of Calculus this immediately implies (\ref{zaq4.5}). $\Box$

%
%
%

\begin{a1}
\label{L6}
Suppose $u$ is a minimizer of $I_{\ep}$ over $W^{2,2}_0(B_1(0))$. There exists $r\eqsim \ep^{\frac{1}{6}}(\log(\ep^{-1}))^{\frac{13}{6}}$ and 
$\xi\in\lt\{1,-1\rt\}$ such that 
\begin{equation}
\label{qza12.7}
\inf\lt\{\xi u(z):z\in B_r(0)\rt\}\geq 1-c\ep^{\frac{1}{6}} (\log(\ep^{-1}))^{\frac{13}{6}}
\end{equation}
\end{a1}

\em Proof. \rm First recall that by Lemma \ref{L0.5}, (\ref{uuzu5}) we know that $I_{\ep}(u)\leq c\ep\log(\ep^{-1})$. 
Let $T_1,T_2,\dots T_N$ be as defined in Lemma \ref{L2}. 
By Lemma \ref{L2} there exists $i\in\lt\{1,2,\dots N\rt\}$ such that 
$T_i$ satisfies (\ref{avv14.7}).

By Lemma \ref{RL2} we know $u\in W^{3,2}(B_{1-2\ep}(0))$. Now by approximation of Sobolev functions (see Theorem 3, section 5.33 \cite{evans2}), for any small $\tau>0$ we can find $\ui\in C^{\infty}(B_{1-2\ep}(0))$ such that  
\begin{equation}
\label{avv50}
\|\ui-u\|_{W^{3,2}(B_{1-2\ep}(0))}<\tau. 
\end{equation}
Since 
\begin{equation}
\label{fineq1}
\int_{B_1(0)} \lt|1-\lt|Du\rt|^2\rt|^2 dx \leq c\ep^2 \log(\ep^{-1}) 
\end{equation}
and 
\begin{equation}
\label{fineqq1}
\int_{B_1(0)} \lt|1-\lt|Du\rt|^2\rt|\lt|D^2 u\rt| dx\leq c\ep\log(\ep^{-1}).
\end{equation}
By Sobolev embedding we have that $u$ is $\frac{1}{2}$-Holder and thus 
\begin{equation}
\label{avv50.5}
\sup\lt\{\lt|u(z)\rt|:z\in \partial B_{1-2\ep}(0)\rt\}\leq c\sqrt{\ep}.
\end{equation}
Now assuming $\tau$ is small enough, as by Sobolev embedding 
$D\ui$ is Holder continuous, $\ui$ must satisfy 
$\sup\lt\{\lt|\ui(z)\rt|:z\in \partial B_{1-2\ep}(0)\rt\}\leq c\sqrt{\ep}$ and  
\begin{eqnarray}
&~&\inf\lt\{\lt|D\ui\lt(z\rt)\rt|:z\in A(0,c\log(\ep^{-1})\ep,1-2\ep)\cap T_i\rt\}>\frac{1}{3}
\text{ and }\nn\\
&~&\qd\qd\sup\lt\{\lt|D\ui\lt(z\rt)\rt|:z\in  A(0,c\log(\ep^{-1})\ep,1-2\ep)\cap T_i\rt\}<3.
\end{eqnarray}
It is also clear that for small enough $\tau$, $\ui$ satisfies $I_{\ep}(\ui)\leq c\ep\log(\ep^{-1})$. \nl

\begin{figure}[h]
\centerline{\psfig{figure=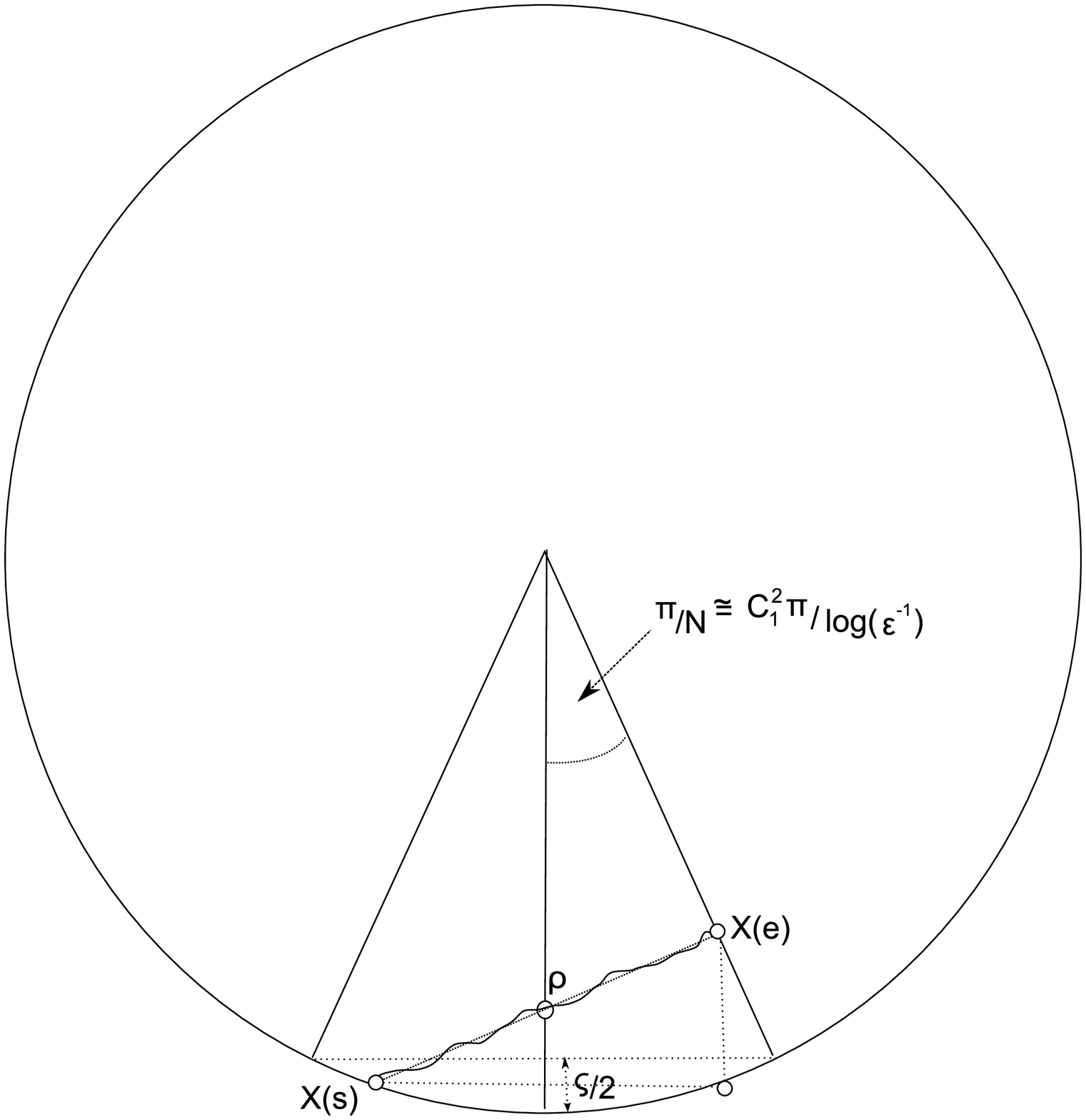,height=9cm,angle=0}}
\caption{}\label{pic}
\end{figure}

\em Step 1. \rm Let $\vthe$ denote the center point of $\partial B_{1-2\ep}(0)\cap T_i$ define 
$\varsigma=2(1-\cos(\frac{\pi}{N}))$, so $\varsigma\eqsim \frac{\cv_1^4\pi^2}{(\log(\ep^{-1}))^2}$. 
Let $\vva=(1-\varsigma)\vthe$. For any set $A$ let $\mathrm{conv}(A)$ denote the 
convex hull of $A$. Note that (see figure \ref{pic}) 
\begin{equation}
\label{nm1}
\mathrm{dist}\lt(\vva,\mathrm{conv}(\partial B_{1-2\ep}(0)\cap T_i)\rt)>\frac{\varsigma}{2}.
\end{equation}

Let
$X:\R\rightarrow \R^2$ be the solution of $X(0)=\vva$ and $\dot{X}(s)=D\ui(X(s))$. Let $\TT_i:=T_i\cap A(0,c\log(\ep^{-1})\ep,1-2\ep)$. 
Let $t_2>0$ be the smallest number such that $X(t_2)\in \partial \TT_i$ and 
let $t_1<0$ be the largest number so that $X(t_1)\in \partial \TT_i$. Let $s\in \lt\{t_1,t_2\rt\}$ be such that 
\begin{equation}
\label{ulazy1}
d(X(s),\partial B_{1-2\ep}(0))=\min\lt\{d(X(t_1)),\partial B_{1-2\ep}(0)),d(X(t_2)),\partial B_{1-2\ep}(0))\rt\}.
\end{equation}
Let $e\in\lt\{t_1,t_2\rt\}\backslash \lt\{s\rt\}$. See figure \ref{pic}.

We will show $X(s)\in \partial B_{1-2\ep}(0)\cap B_{\cv_1^2(\log(\ep^{-1}))^{-1}/2}(\vthe)$ and 
$X(e)\in \partial \TT_i\backslash \partial B_{1-2\ep}(0)$. 

\em Proof of Step 1. \rm  We claim 
\begin{equation}
\label{zaq10}
\cos^{-1}\lt(\frac{X(s)-X(e)}{\lt|X(s)-X(e)\rt|}\cdot \frac{\vthe}{\lt|\vthe\rt|}\rt)\leq \frac{\pi}{2}-\frac{1}{129}.
\end{equation}

Let $\psi=\cos^{-1}\lt(\frac{X(s)-X(e)}{\lt|X(s)-X(e)\rt|}\cdot \frac{\vthe}{\lt|\vthe\rt|}\rt)$. 
Suppose (\ref{zaq10}) not true, i.e.\ $\psi\geq \frac{\pi}{2}-\frac{1}{129}$. Since $X(s),X(e)\not\in B_{\varsigma}(\vthe)$ and by 
(\ref{avv50}), (\ref{fineq1}), (\ref{fineqq1}) $\ui$ satisfies (\ref{qza1}), (\ref{qza2}) for 
$\beta=\ep \log(\ep^{-1})$ so applying Lemma \ref{LL5} we have that by (\ref{zaq4.5})  
\begin{equation}
\label{zaq20}
\vva\in N_{c\ep^{\frac{1}{6}}(\log(\ep^{-1}))^{\frac{7}{6}}}([X(s),X(e)]),
\end{equation}
i.e.\ points $\vva,X(s_2),X(s_1)$ are \em roughly \rm (with error $c\ep^{\frac{1}{6}} (\log(\ep^{-1}))^{\frac{7}{6}}$) aligned, so by (\ref{nm1}) we must have 
$$
X(e)\in \partial \TT_i\backslash \partial B_{1-2\ep}(0)
$$ 
and in particular $\lt|X(e)-X(s)\rt|> \frac{\cv_1^2}{2}(\log(\ep^{-1}))^{-1}$. 
Note also by (\ref{ulazy1}) and by (\ref{zaq20}) we have that 
\begin{equation}
\label{ulazq2}
d(X(s),\partial B_{1-2\ep}(0))\leq c(\log(\ep^{-1}))^{-2}.
\end{equation}
Thus by (\ref{zaq4})
\begin{equation}
\label{qza12}
\lt|\ui(X(e))-\ui(X(s))\rt|\geq \frac{\cv_1^2}{3}(\log(\ep^{-1}))^{-1}.
\end{equation} 

Since $\ui$ is $3$-Lipschitz and $d(X(s),\partial B_{1-2\ep}(0))\leq 2\varsigma$ we have 
$\lt|\ui(X(s))\rt|\leq 6\varsigma\leq \frac{c}{(\log(\ep^{-1}))^2}$. 
Thus by (\ref{qza12}) we have 
\begin{equation}
\label{uzuz1}
\lt|\ui(X(e))\rt|\geq \frac{\cv_1^2}{4} (\log(\ep^{-1}))^{-1}. 
\end{equation}

Now let $L$ be the line parallel to $[X(s),X(e)]$ that passes through $\vva$, by (\ref{zaq4.5}) we can pick 
$\nu\in L\cap B_{\ep^{\frac{1}{6}}(\log(\ep^{-1}))^{\frac{7}{6}}}(X(s))$ 
and let $\mu=(X(e)+\la \vthe\ra)\cap (\nu+\vthe^{\perp})$. Note that by trigonometry 
\begin{equation}
\label{ulazz1}
d(\mu,\partial B_{1-2\ep}(0))\leq d(\nu,\partial B_{1-2\ep}(0))+c(\log(\ep^{-1}))^{-2}.
\end{equation}
And so 
\begin{equation}
\label{ulaz98}
d(\mu,\partial B_{1-2\ep}(0))\leq d(X(s),\partial B_{1-2\ep}(0))+c(\log(\ep^{-1}))^{-2}\overset{(\ref{ulazq2})}{\leq} c(\log(\ep^{-1}))^{-2}.
\end{equation}

Recall we have assumed by contradiction that 
$\psi\geq \frac{\pi}{2}-\frac{1}{129}$. By (\ref{zaq20}) $X(s)$, $\vva$, $X(e)$ are with error $(\ep^{\frac{1}{6}}(\log(\ep^{-1})))^{\frac{7}{6}}$ aligned and by (\ref{ulazy1}) $X(s)$ is closer (or equally close) to $\partial B_{1-2\ep}(0)$ than $X(e)$, so 
$X(s)\cdot \frac{\vthe}{\lt|\vthe\rt|}>X(e)\cdot \frac{\vthe}{\lt|\vthe\rt|}-c\ep^{\frac{1}{6}}(\log(\ep^{-1}))^{\frac{7}{6}}$, hence 
$\psi\leq \frac{\pi}{2}+\frac{1}{129}$. We will denote a triangle with corners at $a,b,c$ by $T(a,b,c)$. Consider the right angle triangle $T(\nu,X(e),\mu)$. Now let $\ti{\psi}$ denote the angle of the corner 
of the triangle $T(\nu,X(e),\mu)$ at $X(e)$. By construction as $\lt|\nu-X(s)\rt|<\ep^{\frac{1}{6}}
(\log(\ep^{-1}))^{\frac{7}{6}}$ so $\lt|\psi-\ti{\psi}\rt|\leq \ep^{\frac{1}{6}}(\log(\ep^{-1}))^{\frac{13}{6}}\leq 
\frac{1}{128}-\frac{1}{129}$, thus $\ti{\psi}\in 
\lt[\frac{\pi}{2}-\frac{1}{128},\frac{\pi}{2}+\frac{1}{128}\rt]$. Thus  
$$
\frac{127}{128}\lt|\nu-X(e)\rt|\leq \lt|\nu-X(e)\rt|\sin(\ti{\psi})\leq \lt|\mu-\nu\rt| \leq 2\pi\cv_1^2 (\log(\ep^{-1}))^{-1}.
$$
So 
\begin{equation}
\label{ulazx3}
\lt|\nu-X(e)\rt|\leq 8\cv_1^2 (\log(\ep^{-1}))^{-1}.
\end{equation}
Thus 
\begin{eqnarray}
\label{ulazz1.5}
\lt|X(e)-\mu\rt|&\leq& \cos(\ti{\psi})\lt|\nu-X(e)\rt|\nn\\
&\overset{(\ref{ulazx3})}{\leq}&8\cv_1^2 (\log(\ep^{-1}))^{-1}\cos\lt(\frac{\pi}{2}-\frac{1}{128}\rt)\nn\\
&\leq& \frac{\cv_1^2(\log(\ep^{-1}))^{-1}}{16}.
\end{eqnarray} 
Hence
\begin{eqnarray}
d(X(e),\partial B_{1-2\ep}(0))&\overset{(\ref{ulazz1.5})}{\leq}&d(\mu,\partial B_{1-2\ep}(0))+\frac{\cv_1^2(\log(\ep^{-1}))^{-1}}{16}\nn\\
&\overset{(\ref{ulaz98})}{\leq}&\frac{\cv_1^2(\log(\ep^{-1}))^{-1}}{16}+c(\log(\ep^{-1}))^{-2}.\nn
\end{eqnarray}
Thus $\lt|\ui(X(e))\rt|\leq \frac{3\cv_1^2(\log(\ep^{-1}))^{-1}}{16}+c\lt(\log(\ep^{-1})\rt)^2$ which is a contradicts 
(\ref{uzuz1}). So (\ref{zaq10}) is established. 

Let $\omega=L\cap (\vthe+\vthe^{\perp})$. 
Consider the right angle triangle $T(\omega,\vva,\vthe)$. By trigonometry we know that 
$\lt|\omega-\vthe\rt|\tan\lt(\frac{\pi}{2}-\psi\rt)=\varsigma$ which 
implies $\lt|\omega-\vthe\rt|\leq 258\varsigma$, hence 
$X(s)\in \partial B_{1-2\ep}(0)\cap B_{\frac{\cv_1^2(\log(\ep^{-1}))^{-1}}{2}}(\vthe)$. As we know already $X(e)\in \partial \TT_i\backslash B_{1-2\ep}(0)$ this 
completes the proof of Step 1. \nl\nl

\em Step 2. \rm We will show  
\begin{equation}
\label{qza27}
\lt|\cos^{-1}\lt(\frac{X(s)}{\lt|X(s)\rt|}\cdot \frac{\lt(X(s)-X(e)\rt)}{\lt|X(s)-X(e)\rt|}\rt)\rt|\leq 
c\ep^{\frac{1}{6}}\log(\ep^{-1})^{\frac{7}{6}}.
\end{equation}

\em Proof of Step 2. \rm Let 
$\theta=\cos^{-1}\lt(\frac{X(s)}{\lt|X(s)\rt|}\cdot \frac{\lt(X(s)-X(e)\rt)}{\lt|X(s)-X(e)\rt|}\rt)$. 
Let 
$$
\kappa=(X(s)+(X(s))^{\perp})\cap \lt(X(e)+\R X(s)\rt). 
$$
Note that the points $X(s), X(e), \kappa$ forms the corners of a right-angle triangle where the angle at the point $X(e)$ is 
$\theta$. Since $\kappa\not\in \TT_i$ and as $\TT_i$ is convex, $\lt[\kappa,X(e)\rt]$ intersects $\partial \TT_i$ at one point only, so let 
$\zeta=(\kappa,X(e))\cap \partial \TT_i$. We claim that $\zeta\in \partial B_{1-2\ep}(0)$. To see this suppose it is not true, 
then the line segment $\lt[\kappa,X(e)\rt]$ must cross one of the flat sides of $\partial \TT_i$. Recall the angle 
at $0$ of the `pie slice' $\TT_i$ is $\frac{2\pi}{N}$. So the angle between $\vthe$ and either of the sides of $\partial \TT_i$ is 
$\frac{\pi}{N}$. However the line segment $[\kappa,X(e)]$ is parallel to the line segment $[0,X(s)]$ so 
$\cos^{-1}\lt(\frac{\vthe}{\lt|\vthe\rt|}\cdot \frac{\kappa-X(e)}{\kappa-X(e)}\rt)< \frac{\pi}{N}$. Now in order for 
$\lt[\kappa,X(e)\rt]$ to cross the flat sides of $\partial \TT_i$ without first intersecting $\partial B_{1-2\ep}(0)$ it 
has to make a larger angle with $\vthe$ than the flat sides of $\partial \TT_i$ so this a contradiction. 
Thus the claim is established and we have $\cos(\theta)\lt|X(s)-X(e)\rt|\geq \lt|X(e)-\zeta\rt|$.

Now since $X(s)\in \partial B_{1-2\ep}(0)$ so $\lt|\ui(X(s))\rt|\leq c\sqrt{\ep}$ and thus 
\begin{eqnarray}
\label{zeq30}
\ui(X(e))&\overset{(\ref{zaq4})}{\geq}& (1-c(\ep\log(\ep^{-1}))^{\frac{1}{3}})\lt|X(e)-X(s)\rt|-c(\ep\log(\ep^{-1}))^{\frac{1}{3}}\nn\\
&\geq& \frac{\lt|X(e)-\zeta\rt|}{\cos \theta}-c(\ep\log(\ep^{-1}))^{\frac{1}{3}}.
\end{eqnarray}

By Lemma \ref{L5} there exists a line segment $\Gamma\subset \TT_i$ parallel to $[X(e),\zeta]$ whose end points are within 
$(\ep \log(\ep^{-1}))^{\frac{1}{3}}$ of $X(e),\zeta$ and for which 
$\sup\lt\{\lt|\lt|D \ui(z)\rt|-1\rt|:z\in \Gamma\rt\}\leq c(\ep\log(\ep^{-1}))^{\frac{1}{3}}$. Let $a,b$ be the end points of 
$\Gamma$, so by the fundamental theorem of Calculus, $\lt|\ui(a)-\ui(b)\rt|\leq (1+ c(\ep\log(\ep^{-1}))^{\frac{1}{3}})
\lt|a-b\rt|$. Since $\ui$ is Lipschitz on $\TT_i$ and $\lt|\ui(\zeta)\rt|\leq c\sqrt{\ep}$ we have that  $\lt|\ui(X(e))\rt|\leq (1+ c(\ep\log(\ep^{-1}))^{\frac{1}{3}})\lt|X(e)-\zeta\rt|$, thus putting this together with (\ref{zeq30}) we have 
\begin{equation}
\label{zeq31}
\lt|X(e)-\zeta\rt|\geq \frac{\lt|X(e)-\zeta\rt|}{(1+ c(\ep\log(\ep^{-1}))^{\frac{1}{3}}) \cos \theta}-c(\ep\log(\ep^{-1}))^{\frac{1}{3}}.
\end{equation}
Recall $B_{\varsigma}(\vva)\subset \TT_i$ and as we know $X(s)$ is closer to $\partial B_{1-2\ep}(0)$ than 
$X(e)$, so by (\ref{zaq20}) we have that $\lt|X(e)-\zeta\rt|\geq 
\frac{\varsigma}{2}$, so 
by (\ref{zeq31}) we have $\cos(\theta)\geq 1-c\ep^{\frac{1}{3}}(\log(\ep^{-1}))^{\frac{7}{3}}$ which implies 
$\lt|\theta\rt|\leq c\ep^{\frac{1}{6}}(\log(\ep^{-1}))^{\frac{7}{6}}$ and this completes the proof of Step 2. $\Box$ \nl

\em Proof of Lemma completed. \rm By Step 1 we know $X(s)\in B_{\frac{\cv_1^2 (\log(\ep^{-1}))^{-1}}{2}}(\vthe)$, so the angle between 
the line segment $[X(s),0]$ and the sides of $\partial \TT_i$ is at least $\cv_1^2(\log(\ep^{-1})^{-1}/4$. So if we 
consider the triangle $T(0,X(s),X(e))$. Let $\eta$ be the angle of the triangle at corner $0$, so 
$\eta\geq \frac{\cv_1^2(\log(\ep^{-1}))^{-1}}{4}$. Recall the angle at corner $X(s)$ is $\theta$ and by (\ref{qza27}) 
$\theta\leq c\ep^{\frac{1}{6}}(\log(\ep^{-1}))^{\frac{7}{6}}$. So by the law of sins, $\frac{\lt|X(e)\rt|}{\sin\theta}
=\frac{\lt|X(e)-X(s)\rt|}{\sin\eta}$. So 
\begin{equation}
\label{uzuu2}
\lt|X(e)\rt|\leq \frac{2\sin \theta}{\sin\eta}\leq 
c\ep^{\frac{1}{6}}(\log(\ep^{-1}))^{\frac{13}{6}}. 
\end{equation}
Now as noted previously, (\ref{avv50.5}) and (\ref{avv50}), $\lt|\ui(X(s))\rt|\leq c\sqrt{\ep}$. So by (\ref{zaq4}) we have that 
\begin{eqnarray}
\label{ulzy12}
\lt|\ui(X(e))\rt|&\geq& (1-(\ep\log(\ep^{-1})^{\frac{1}{3}})\lt|X(e)-X(s)\rt|
-c(\ep\log(\ep^{-1}))^{\frac{1}{3}}\nn\\
&\geq& 
(1-(\ep\log(\ep^{-1}))^{\frac{1}{3}})d(X(e),\partial B_{1-2\ep}(0))-c(\ep\log(\ep^{-1}))^{\frac{1}{3}}\nn\\
&\geq& 1-c\ep^{\frac{1}{6}}(\log(\ep^{-1}))^{\frac{13}{6}}.
\end{eqnarray}
So we must have $r\in (\lt|X(e)\rt|+\frac{1}{2}\ep^{\frac{1}{6}}(\log(\ep^{-1}))^{\frac{13}{6}},
\lt|X(e)\rt|+c\ep^{\frac{1}{6}}(\log(\ep^{-1}))^{\frac{13}{6}})$ such that 
$$
\int_{\partial B_r(0)} \lt|1-\lt|D\ui\rt|^2\rt| dH^1 z\overset{(\ref{fineq1}),(\ref{avv50})}{\leq} c\ep^{\frac{5}{6}}(\log(\ep^{-1}))^{-\frac{10}{6}}.
$$ 
By the fundamental theorem of Calculus was have that 
\begin{equation}
\label{qw1}
\lt|\ui(x)-\ui(y)\rt|\leq c\ep^{\frac{5}{6}}(\log(\ep^{-1}))^{-\frac{10}{6}}\text{ for all }
x,y\in \partial B_r(0).
\end{equation}

Let $\xi=\frac{\ui(X(e))}{\lt|\ui(X(e))\rt|}$. Pick $z\in \partial B_r(0)\cap \TT_i$, since $\ui$ is Lipschitz on $\TT_i$ we know 
\begin{equation}
\label{uzu1}
\lt|\ui(z)-\ui(X(e))\rt|\leq  c\ep^{\frac{1}{6}}(\log(\ep^{-1}))^{\frac{13}{6}}. 
\end{equation}
Thus for any $x\in \partial B_r(0)$ 
\begin{equation}
\label{ulzq1}
\xi \ui(x)\overset{(\ref{uzu1})(\ref{qw1})}{\geq} \xi \ui(X(e))-c\ep^{\frac{1}{6}}(\log(\ep^{-1}))^{\frac{13}{6}}
\overset{(\ref{ulzy12})}{\geq} 1-c\ep^{\frac{1}{6}}(\log(\ep^{-1}))^{\frac{13}{6}},
\end{equation}
together with (\ref{avv50}) (using the fact that (\ref{avv50}) implies $\|\ui-u\|_{L^{\infty}(B_{1-2\ep}(0))}\leq c\ep$) 
this completes the proof of Lemma \ref{L6}.\nl

\em Proof of Theorem  completed. \rm Let $r\eqsim \ep^{\frac{1}{6}}(\log(\ep^{-1}))^{\frac{13}{6}}$, $\xi\in\lt\{-1,1\rt\}$ be the numbers that 
satisfy (\ref{qza12.7}) from Lemma \ref{L6}. Let $A(x)=\frac{x}{\lt|x\rt|}$ note $\lt|DA(x)\rt|\leq \frac{c}{\lt|x\rt|}$. 
Note by Fubini 
\begin{eqnarray}
&~&\int_{B_r(0)}\int_{B_1(0)} \lt|1-\lt|Du(z)\rt|^2\rt|\lt|DA(x-z)\rt| dz dx\nn\\
&~&\qd\qd\qd=\int_{B_1(0)}\lt(\int_{B_r(0)} \lt|DA(x-z)\rt| dx\rt) \lt(1-\lt|Du(z)\rt|^2\rt) dz\nn\\
&~&\qd\qd\qd\leq c\ep\sqrt{\log(\ep^{-1})}.
\end{eqnarray}
So there must exist a set $G\subset B_r(0)$ with $\lt|G\rt|\geq \ep^{\frac{1}{3}}(\log(\ep^{-1}))^{\frac{13}{3}}$ such 
that if $x\in G$ we have 
\begin{equation}
\int_{B_1(0)} \lt|1-\lt|Du(z)\rt|^2\rt|\lt|DA(x-z)\rt| dz\leq c\ep^{\frac{1}{3}}.
\end{equation}
For $\theta\in S^1$, $y\in \R^2$ define $l_{\theta}^y:=y+\R_{+}\theta$. Pick $x\in G$, by the Co-area formula 
$$
\int_{\psi\in S^1} \int_{l_{\psi}^{x}} \lt|1-\lt|Du(z)\rt|^2\rt| dH^1 z dH^1 \psi\leq 
c\ep^{\frac{1}{3}}.
$$
For each $\psi\in S^1$ let $x_{\psi}=\partial B_r(0)\cap l_{\psi}^x$, $y_{\psi}=\partial B_1(0)\cap l_{\psi}^x$ and 
$e_{\psi}=\int_{l^x_{\psi}} \lt|1-\lt|Du(z)\rt|^2\rt| dH^1 z$. So 
\begin{eqnarray}
\label{ulzq2}
\int_{[x_{\psi},y_{\psi}]} \lt|Du(z)+\xi\psi\rt|^2 dH^1 z &=&
\int_{[x_{\psi},y_{\psi}]} \lt|Du(z)\rt|^2+2\xi Du(z)\cdot\psi+1 dH^1 z \nn\\
&\leq&2\lt|y_{\psi}-x_{\psi}\rt|-2 \xi u(x_{\psi})+ c e_{\psi}\nn\\
&\overset{(\ref{qza12.7})}{\leq}&c\ep^{\frac{1}{6}}(\log(\ep^{-1}))^{\frac{13}{6}}+ c e_{\psi}.
\end{eqnarray}

Thus 
\begin{eqnarray}
\int_{B_1(0)\backslash B_r(x)} \lt|Du(z)+\xi \frac{z}{\lt|z\rt|}\rt|^2 dz&\leq&
\int_{B_1(0)\backslash B_r(x)} \lt|Du(z)+\xi \frac{z}{\lt|z\rt|}\rt|^2\lt|DA(x-z)\rt| dz\nn\\
&\leq&\int_{S^1} \int_{[x_{\psi},y_{\psi}]} \lt|Du(z)+\xi \psi\rt| dH^1 z dH^1 \psi\nn\\
&\overset{(\ref{ulzq2})}{\leq}&c\ep^{\frac{1}{6}}(\log(\ep^{-1}))^{\frac{13}{6}}+c\int_{S^1} e_{\psi} dH^1 \psi\nn\\
&\leq&c\ep^{\frac{1}{6}}(\log(\ep^{-1}))^{\frac{13}{6}}.\nn
\end{eqnarray}

Hence 
\begin{eqnarray}
\int_{B_1(0)} \lt|Du(z)+\xi \frac{z}{\lt|z\rt|}\rt|^2 dz&\leq&
\int_{B_r(0)} \lt|Du(z)+\xi \frac{z}{\lt|z\rt|}\rt|^2 dz+c\ep^{\frac{1}{6}}(\log(\ep^{-1}))^{\frac{13}{6}}\nn\\
&\leq&c\int_{B_r(0)} \lt|1-\lt|\lt|Du(z)\rt|-1\rt|\rt|^2 dz+c\ep^{\frac{1}{6}}(\log(\ep^{-1}))^{\frac{13}{6}} \nn\\
&\leq&c\ep^{\frac{1}{6}}(\log(\ep^{-1}))^{\frac{13}{6}}.\nn \;\;\;\;\;\;\Box
\end{eqnarray}

\bf Acknowledgments. \rm I would like to thank Michael Goldberg for helpful discussions and the anonymous referee 
for careful reading, several good suggestions and indicating the possibility of a simplification of part of Lemma 2.

\end{document}